\documentclass[12pt]{amsart}
\usepackage{amssymb}
\usepackage[mathscr]{eucal}

\usepackage{colordvi}
\usepackage{graphicx}
\usepackage{ulem}    
\usepackage{rotating}
\usepackage{multirow}

\usepackage{icomma}
\usepackage[ps2pdf]{hyperref}
\usepackage{url}

\usepackage{vmargin}
\setpapersize{USletter}
\setmargrb{1in}{1in}{1in}{1in} 
 
\numberwithin{equation}{section}
\newtheorem{theorem}{Theorem}[section]
\newtheorem{defn}[theorem]{Definition}
\newtheorem{proposition}[theorem]{Proposition}

\newtheorem{conj}[theorem]{Conjecture}
\newtheorem{Shapiroconj}[theorem]{Shapiro Conjecture}
\newtheorem{sconj}[theorem]{Secant Conjecture}
\newtheorem{monconj}[theorem]{Monotone Conjecture}
\newtheorem{monsecconj}[theorem]{Monotone Secant Conjecture}
\newtheorem{example}[theorem]{Example}
\newtheorem{remark}[theorem]{Remark}

\theoremstyle{remark}

\newcounter{FNC}[page]
\def\fauxfootnote#1{{\addtocounter{FNC}{2}\Magenta{$^\fnsymbol{FNC}$}%
     \let\thefootnote\relax\footnotetext{\Magenta{$^\fnsymbol{FNC}$#1}}}}

\newcommand{\adot}{a_{\bullet}}
\newcommand{\Fdot}{F_{\bullet}}
\newcommand{\Edot}{E_{\bullet}}
\newcommand{\Fl}{\mathbb{F}\ell}
\newcommand{\Gr}{{\rm Gr}}

\newcommand{\C}{{\mathbb{C}}}

\renewcommand{\P}{{\mathbb{P}}}

\newcommand{\R}{{\mathbb{R}}}

\newcommand{\W}{\Purple{2}}
\newcommand{\X}{\Blue{3}}

\newcommand{\bsigma}{{\boldsymbol\sigma}}

\newcommand{\Span}{\mbox{\rm span}}

\newcommand{\DeCo}[1]{\Blue{#1}}
\newcommand{\demph}[1]{\DeCo{{\sl #1}}}

\def\rgbColor#1#2{#2}
\newcommand{\M}[1]{\rgbColor{0.8 0 0.8}{#1}}
\newcommand{\B}[1]{\Blue{#1}}
\newcommand{\G}[1]{\rgbColor{0 0.5 0}{#1}}

\title[Monotone Secant Conjecture]{The Monotone Secant Conjecture in the\\ real Schubert calculus}

\author[Hein]{Nickolas Hein}
\address{Nickolas Hein \\
         Department of Mathematics\\
         University of Nebraska at Kearney\\
         Kearney\\
         Nebraska \ 68849\\
         USA}
\email{heinnj@unk.edu}
\urladdr{\url{http://www.unk.edu/academics/math/faculty/About_Nickolas_Hein/}}
\author[Hillar]{Christopher J. Hillar}
\address{Christopher J. Hillar \\
         Mathematical Sciences Research Institute\\
         17 Gauss Way\\
         Berkeley, CA 94720-5070\\
         USA}
\email{chillar@msri.org}
\urladdr{\url{http://www.msri.org/people/members/chillar}}
\author[Mart\'in del Campo]{Abraham Mart\'in del Campo}
\address{Abraham Mart\'in del Campo\\
         IST Austria\\
         Am Campus 1\\
         3400 Klosterneuburg\\
         Austria
         }
\email{abraham.mc@ist.ac.at}
\urladdr{\url{http://pub.ist.ac.at/~adelcampo/}}
\author[Sottile]{Frank Sottile}
\address{Frank Sottile \\
         Department of Mathematics\\
         Texas A\&M University\\
         College Station\\
         Texas \ 77843\\
         USA}
\email{sottile@math.tamu.edu}
\urladdr{\url{http://www.math.tamu.edu/~sottile}}
\author[Teitler]{Zach Teitler}
\address{Zach Teitler \\
         Department of Mathematics\\
         Boise State University\\
         Boise\\
         Idaho \ 83725\\
         USA}
\email{zteitler@math.boisestate.edu}
\urladdr{\url{http://math.boisestate.edu/~zteitler}}

\thanks{Research supported in part by NSF grants DMS-0701050, DMS-0915211, and DMS-1001615.}
\thanks{Research of Hillar  supported in part by an NSF Postdoctoral Fellowship 
        and an NSA Young Investigator grant}
\thanks{This research conducted in part on computers provided by NSF SCREMS grant DMS-0922866}

\begin{document}

\begin{abstract}
 The monotone secant conjecture posits a rich class of polynomial systems, all of whose
 solutions are real.
 These systems come from the Schubert calculus on flag manifolds, and the monotone secant
 conjecture is a compelling generalization of the Shapiro conjecture for
 Grassmannians (Theorem of Mukhin, Tarasov, and Varchenko).
 We present some theoretical evidence for this conjecture, as well as computational evidence
 obtained by 1.9 teraHertz-years of computing, and we discuss some 
 of the phenomena we observed in our data.
\end{abstract}

\maketitle
%
\section{Introduction}

A system of real polynomial equations with finitely many solutions has some, but likely
not all, of its solutions real. 
In fact, sometimes the structure of the equations implies an upper bound on the number of real
solutions~\cite{BBS, Kh91}, ensuring that not all solutions are real. 
The monotone secant conjecture posits a family of systems of
polynomial equations with the extreme property that all of their solutions are real. 

The Shapiro conjecture asserts that a zero-dimensional intersection of Schubert
subvarieties of a Grassmannian consists only of real points provided that the Schubert
varieties are given by flags tangent to a real rational normal curve. 
While the statement concerns the Schubert calculus on
Grassmannians, its proofs involve complex
analysis~\cite{EG_02,EG11} or integrable systems and representation
theory~\cite{MTV_Annals,MTV_JAMS}. 
A complete story of this conjecture and its proof can be found in~\cite{So_FRSC}.

The Shapiro conjecture is false for non-Grassmannian flag manifolds, but in a
very interesting manner. 
This failure was first noticed in~\cite{So_Shap} and systematic computer experimentation
suggested a correction, the monotone conjecture~\cite{RSSS, So_fulton},
that appears to be valid for flag manifolds of type A. 
Eremenko, Gabrielov, Shapiro, and
Vainshtein~\cite{EGSV} proved a result that implies the monotone conjecture for some
manifolds of two-step flags and concerns codimension-two subspaces that meet flags
which are secant to the rational normal curve along disjoint intervals.
This suggested the secant conjecture, which
asserts that an intersection of Schubert varieties in a Grassmannian is transverse with all
points real, provided that the Schubert varieties are defined by flags secant to a
rational normal curve along disjoint intervals.
This was posed and evidence was presented for its validity in~\cite{FRSC_Sec}.

The monotone secant conjecture is a common extension of both the
monotone conjecture and the secant conjecture.
It is also the last of the conjectures our group has made concerning reality in Schubert
calculus of osculating flags.
In addition to those mentioned, there is a version of the Shapiro conjecture for the orthogonal
Grassmannian which was proven by Purbhoo~\cite{Purbhoo}, and a version for the Lagrangian
Grassmannian described in~\cite[Ch.~14.2]{IHP}.
Exploratory computations in other flag manifolds suggest there is no regularity in the
number of real solutions to Schubert calculus problems given by osculating flags.

We give here an open instance of the monotone secant conjecture, expressed as a system of
polynomial equations in local coordinates for the variety of flags $E_2\subset E_3$ in
$\C^5$, where $\dim E_i=i$. Let $x_1,\ldots,x_8$ be indeterminates and consider the
polynomials 
 \begin{equation}\label{Eq:F235}
   f(s,t,u;x)\ :=\ \det\left(\begin{matrix}
                     1&0&x_1&x_2&x_3\\ 
                     0&1&x_4&x_5&x_6\\\hline
                     1&s&s^2&s^3&s^4\rule{0pt}{12pt}\\
                     1&t&t^2&t^3&t^4\\
                     1&u&u^2&u^3&u^4
                   \end{matrix}\right)\ , 
                   \hspace{10pt}
  g(v,w;x)\ :=\ \det\left(\begin{matrix}
                     1&0&x_1&x_2&x_3\\ 
                     0&1&x_4&x_5&x_6\\
                     0&0&1&x_7&x_8\\\hline
                     1&v&v^2&v^3&v^4\rule{0pt}{12pt}\\
                     1&w&w^2&w^3&w^4\\
		    \end{matrix}\right),
\end{equation}
which depend upon parameters $s,t,u$ and $v,w$ respectively.

\begin{conj}\label{C:first}
  Let $\Blue{s_1<t_1<u_1<\dotsb< s_4< t_4< u_4}\,
    <\, \G{v_1<w_1<\dotsb< v_4\,<\, w_4}$ be real numbers.
  Then the system of polynomial equations
 \begin{eqnarray}\label{Eq:polysys}
   \B{f(s_1,t_1,u_1; x)\ = \ f(s_2,t_2,u_2; x)\ =\ f(s_3,t_3,u_3; x)\ =\ f(s_4,t_4,u_4; x)} &=& 0  \\
   \nonumber\G{g(v_1,w_1; x) \ = \ g(v_2,w_2; x) \ = \ g(v_3,w_3; x) \ = \ g(v_4,w_4; x)} &=& 0
 \end{eqnarray}
  has twelve solutions, and all of them are real.
\end{conj}

These equations have geometric meaning.
Let $E_2$ be the span of the first two rows of either matrix and $E_3$ the span of the
first three rows of the second matrix that defines $g$.
Then $E_2\subset E_3$ is a general flag.
If we let $F_3(s,t,u)$ be the span of the last three rows of the matrix for $f$, then this
is a 3-plane that is secant to the rational normal curve $\gamma\colon y\mapsto
(1,y,y^2,y^3,y^4)$ at the points $\gamma(s),\gamma(t),\gamma(u)$.
The equation $f(s,t,u;x)=0$ is the condition that $E_2$ meet $F_3(s,t,u)$ non-trivially.
Smilarly, if $F_2(u,v)$ is the span of the last two rows of the second matrix, which is 
2-plane secant to $\gamma$ at $\gamma(u)$ and $\gamma(v)$, then the equation  $g(v,w;x)=0$
is the condition that $E_3$ meet $F_2(u,v)$ non-trivially.

The monotonicity hypothesis is that the four 3-planes given by $s_i,t_i, u_i$ are
secant along intervals $\B{[s_i,u_i]}$ which are pairwise disjoint and 
occur before the pairwise disjoint intervals $\G{[v_i, w_i]}$ where the 2-planes
are secant. 
If the order of the intervals $\B{[s_4,u_4]}$ and $\G{[v_1, w_1]}$ is reversed, the
evaluation is no longer monotone. 
We tested $3,000,000$ instances of Conjecture~\ref{C:first}, finding only real solutions.
In contrast, we tested $21,000,000$ with the monotonicity condition violated, of which
$18,085,537$ had some non real solutions.

We formulate the monotone secant conjecture, explain its relation to the other reality 
conjectures, describe data supporting it from a large computational experiment, and discuss some
features observed in our data that go beyond the monotone secant conjecture.
This experiment verified the monotone secant conjecture in each of the 768,846,000
instances it tested.
We have created a website~\cite{FRSC} for viewing the data online.
This includes pages for browsing the data and viewing the results for each Schubert
problem.
We only sketch the other reality conjectures, for they are described in the cited literature,
and we also only sketch the design and execution of this experiment, for 
the purpose of the paper~\cite{Exp-FRSC} was to present the software environment we 
developed for such distributed computational experiments.

This paper is organized as follows. 
In Section~\ref{Sec:fourlines} we illustrate the main point of the monotone secant
conjecture through the classical problem of four lines.  
Section \ref{Sec:background} is a primer on flag manifolds and contains a precise
statement of the monotone secant conjecture while also explaining its relation to the
Shapiro, secant, and monotone conjectures.
In Section \ref{Sec:results} we expand on the relation between the monotone secant and monotone
conjectures, discuss the experimental evidence for the monotone secant conjecture, and some 
phenomena we observed in our data.
Lastly, in Section \ref{Sec:method} we sketch the methods we used to test the conjecture. 

%
\section{The problem of four lines}\label{Sec:fourlines} 
The classical problem of four lines asks for the
finitely many lines $\B{m}$ that meet four given general lines $\G{\ell_1},\, \M{\ell_2},\,
\Red{\ell_3},\, \ell_4$ in (projective) three-space. 
This has a pleasing synthetic solution, which leads to the first interesting case of the
monotone secant conjecture.

Three general lines 
$\G{\ell_1},\,\M{\ell_2},\,\Red{\ell_3}$ lie in one ruling of a doubly-ruled
quadric surface $\Brown{Q}$, with the other ruling consisting of all lines that meet the
first three. 
The line $\ell_4$ meets $\Brown{Q}$ in two points, and through each of these points there
is a line in the second ruling. 
These two lines, $\B{m_1}$ and $\B{m_2}$, are the solutions to this problem.
If the lines $\G{\ell_1},\, \M{\ell_2},\, \Red{\ell_3},\, \ell_4$ are real,
then so is \Brown{$Q$}, but the intersection of \Brown{$Q$} with $\ell_4$ is 
either two real points or a pair of complex conjugate points. 
In the first case, the problem of four lines has two real solutions, while in the second, 
it has no real solutions.  

Let us consider a variant in the manifold of flags consisting of a line $m$ lying on a plane
$M$ in 3-space, $m\subset M$.
Consider the Schubert problem in which $m$ meets three lines 
$\G{\ell_1},\,\M{\ell_2},\,\Red{\ell_3}$ and $M$ contains two
points, $p$ and $q$.
Then $M$ contains the affine span $\overline{p,q}$ of $p$ and $q$.
Since $m\subset M$, we must have that $m$ also meets $\overline{p,q}$ and is therefore a
solution to the problem of four lines given by 
$\G{\ell_1},\,\M{\ell_2},\,\Red{\ell_3}$ and $\overline{p,q}$.
As $M$ is spanned by $m$ and $\overline{p,q}$, we see that solving this auxiliary problem of
four lines solves our original Schubert problem.
Furthermore, if the lines and points are real, then a solution $m\subset M$ is real if and only
if $m$ is real.

Suppose that the three lines are secant to a rational normal curve $\gamma$
along disjoint intervals and the 
points are $p=\gamma(s)$ and $q=\gamma(t)$, 
which do not lie in 
any interval of secancy.
There are two possible combinatorial placements of the two points.
Removing the three intervals of secancy from $\gamma$ results in three disjoint intervals along
$\gamma\simeq\R\P^1$.
Either both points $\gamma(s)$ and $\gamma(t)$ lie in the same interval or they lie in
different intervals.
We examine each case.

\begin{figure}[htb]
\[
  \begin{picture}(320,110)(0,7)
   \put(0,7){\includegraphics[height=119pt,viewport=5 78 430 240,clip]{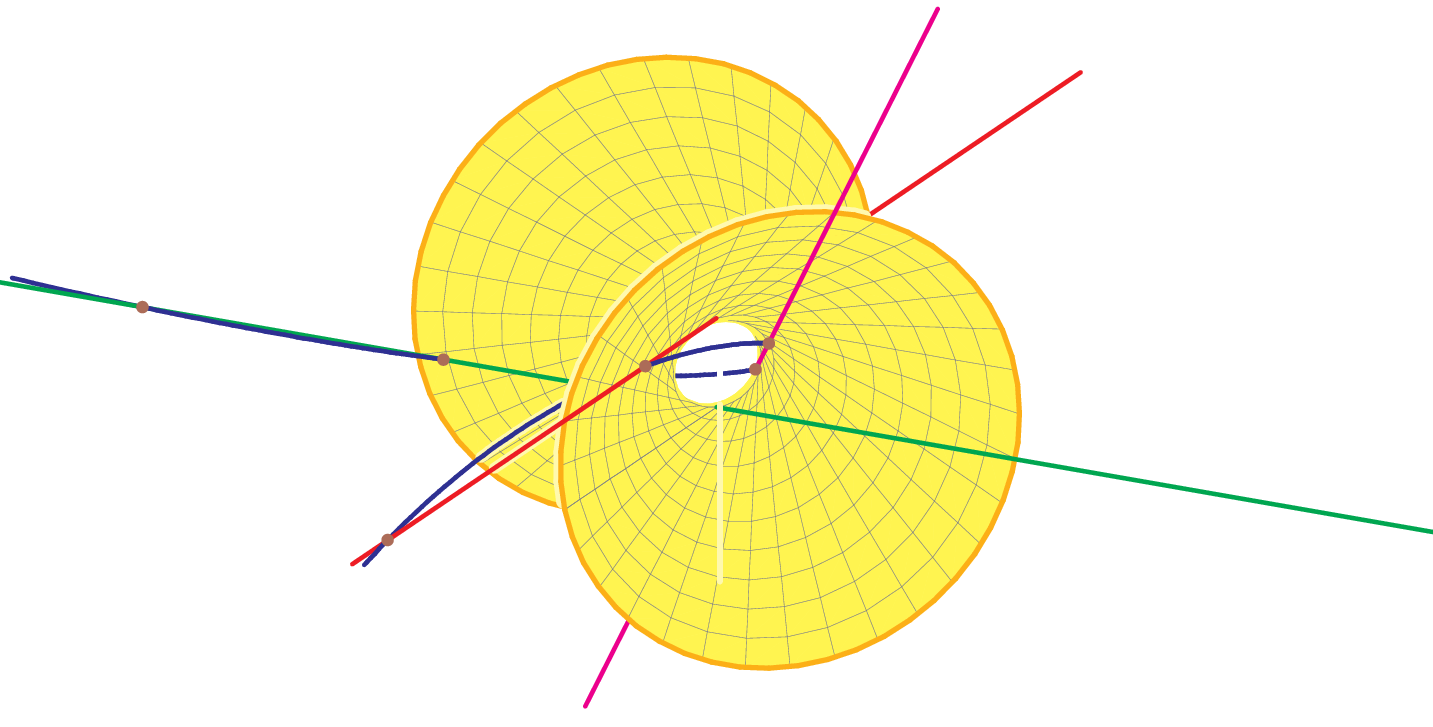}}
   \thicklines
   \put(119,49){$\gamma$}\put(128,51){\vector(2,-1){25}}
   \put(250,0){\vector(0,1){82.3}}
   \put(246,-10){$I$}
   \put(320,90){\Brown{$Q$}}
  \end{picture}
\]
\caption{The problem of four secant lines.}\vspace{-10pt}
\label{F:secant_four}
\end{figure}
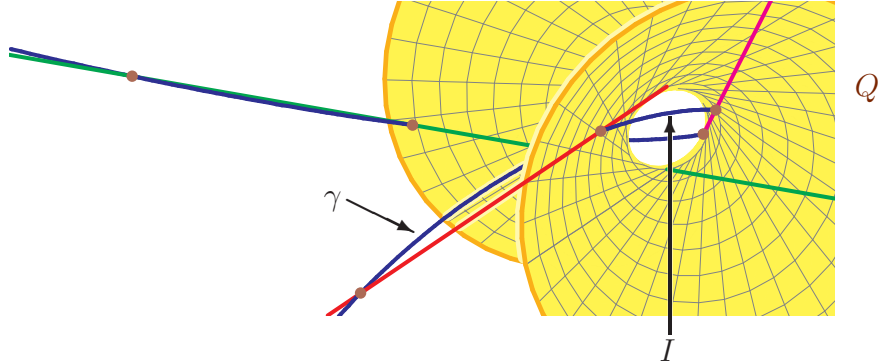

\begin{figure}[htb]
\[
  \begin{picture}(310,110)(0,15)
\put(0,7){\includegraphics[height=119pt,viewport=5 78 430 240,clip]{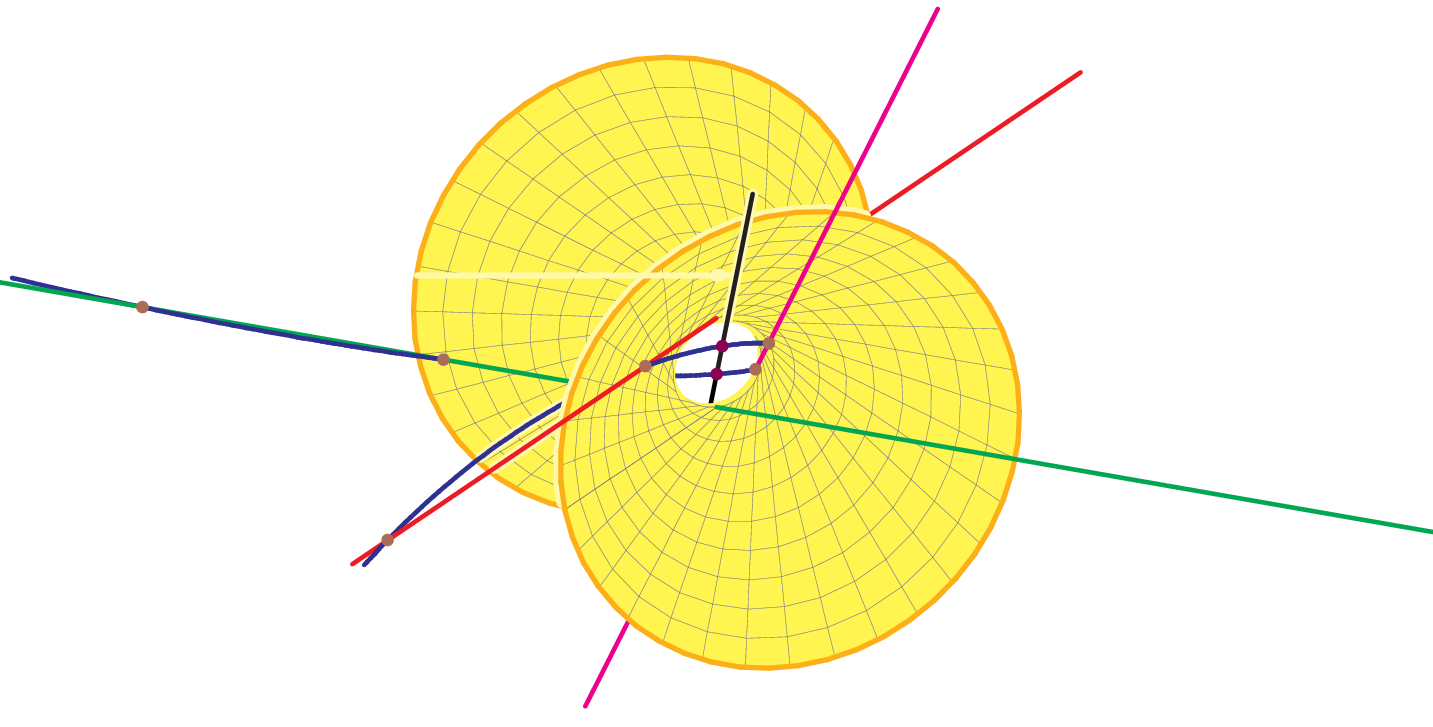}}
   \thicklines
   \put(119,49){$\gamma$}\put(128,51){\vector(2,-1){25}}
   \put(119,109){\vector(1,0){134}}
   \put(90,106){$\ell(s,t)$}
   \put(320,90){\Brown{$Q$}}
  \end{picture}
\]
\caption{A non-monotone evaluation.}
\label{F:not_monotone}
\end{figure}
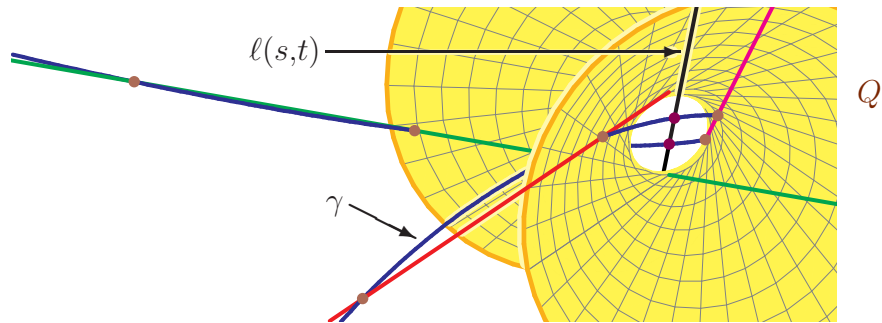

Fixing 
secant lines $\G{\ell_1},\, \M{\ell_2},\, \Red{\ell_3}$, the quadric
\Brown{$Q$} described above is a 
hyperboloid of one sheet.
This is displayed in Figures~\ref{F:secant_four} and~\ref{F:not_monotone}, along with $\gamma$
and the lines.
Suppose that $\gamma(s)$ and $\gamma(t)$ lie in the same interval, say $I$, as indicated in
Figure~\ref{F:secant_four}.
Then the secant line they span, $\ell(s,t)$, lies in the direction of $I$ and meets the
hyperboloid \Brown{$Q$} in two real points.
Thus, in this first case, our Schubert problem has two real solutions.
(This is also an instance of the secant conjecture, which holds for this problem of four
lines~\cite[\S4]{FRSC_Sec}.)

In the second case, where the points $\gamma(s)$ and $\gamma(t)$ do not lie in the same
interval, it is possible to have no real solutions.
Consider the choice of points  $\gamma(s)$ and $\gamma(t)$ as shown in 
Figure~\ref{F:not_monotone}, so that the secant line $\ell(s,t)$
does not meet the quadric \Brown{$Q$}.
By our previous analysis, there will be no real lines $m$ meeting these four secant
lines, and therefore no real solutions $m\subset M$ to our Schubert problem.

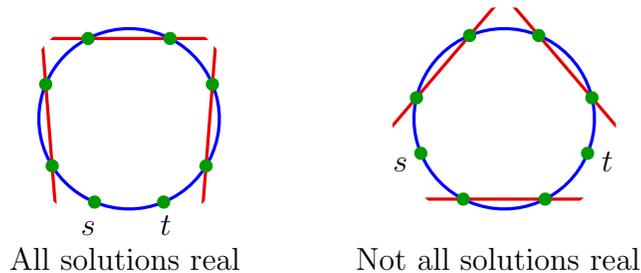
\begin{figure}[htb]
\[
   \begin{picture}(88,90)(-8,-17)
    \put(0,0){\includegraphics{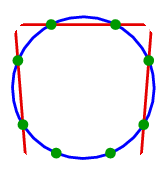}}
    \put(22,-4){$s$} \put(52,-4){$t$}
    \put(-5,-17){All solutions real}
   \end{picture}
    \qquad\qquad
   \begin{picture}(96,90)(-15,-17)
    \put(0,0){\includegraphics{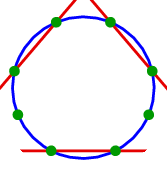}}
    \put(-2,20){$s$} \put(77,20){$t$}
    \put(-16,-17){Not all solutions real}
   \end{picture}
\]
\caption{Schematic for the secancies.}
\label{F:schematic}
\end{figure}

We conclude that the positions of the points $\gamma(s),\gamma(t)$ relative to the other
intervals of secancy may affect whether or not the solutions are real.
The schematic in Figure~\ref{F:schematic}
illustrates the relative positions of the secancies along $\gamma$ (which is homeomorphic
to the circle).
The idea behind the monotone secant conjecture is to attach to each interval the dimension
of that part of the flag which it affects.
This is 1 for $m$ and 2 for $M$.
The schematic on the left has labels $1,1,1,2,2$, reading clockwise, starting just past
the point $s$, while that on the right reads $1,1,2,1,2$.
The labels increase monotonically in the first and do not in the second.

%
\section{Background}\label{Sec:background}

We develop the background for the statement of the monotone secant conjecture, defining
flag varieties and their Schubert problems. 
More may be found in the book of Fulton~\cite{Fu97}.
Fix positive integers $\DeCo{\adot}:=(a_1 < \cdots < a_k)$ and $n$ with $a_k<n$.
A \demph{flag $\Edot$ of type $\adot$} is a sequence of subspaces
\[
   \Edot\ \colon\ \{0\}\subset E_{a_1}\ \subset\ E_{a_2}\ \subset\ 
      \dotsb\ \subset\ E_{a_k}\ \subset\ \C^n\,,
   \qquad\mbox{where\ }\dim(E_{a_i})=a_i\,.
\]
The set of all such flags forms the 
\demph{flag manifold $\Fl(\adot;n)$}, which has dimension 
$\DeCo{\dim(\adot)}:=\sum_{i=1}^k(n-a_i)(a_i-a_{i-1})$, where $a_0:=0$.
When $\adot=(a)$ is a singleton, $\Fl(\adot;n)$ is the \demph{Grassmannian} 
of $a$-planes in $\C^n$, written \DeCo{$\Gr(a,n)$}. 
Flags of type $1<2<\cdots<n-1$ in $\C^n$ are \demph{complete}. 

The positions of flags relative to a fixed complete flag $\Fdot$ stratify
$\Fl(\adot;n)$ into topological cells whose closures are \demph{Schubert varieties}. 
These positions are indexed by certain permutations. 
The \demph{descent set $\delta(\sigma)$} of a permutation $\sigma\in S_n$ is the set of
numbers $i$ such that $\sigma(i)>\sigma(i{+}1)$.
For a permutation $\sigma\in S_n$ with $\delta(\sigma)\subset \adot$,
the Schubert variety $X_\sigma \Fdot$ is
 \[
   X_\sigma \Fdot \ =\ 
   \{ \Edot\in\Fl(\adot;n) \mid 
   \dim E_{a_i}\cap F_j \geq \#\{ l\leq i \mid j+\sigma(l)>n\}\ \ \forall i,j \}.
 \]
Flags $\Edot$ in $X_\sigma \Fdot$ have position $\sigma$ relative to $\Fdot$.
A permutation $\sigma$ with descent set contained in $\adot$ is a \demph{Schubert
  condition} on flags of type $\adot$.  
The Schubert variety $X_\sigma \Fdot$ is irreducible with codimension 
 $\ell(\sigma):= |\{i<j \mid \sigma(i)>\sigma(j)\}|$. 
A \demph{Schubert problem } for $\Fl(\adot;n)$ is a list 
$\DeCo{\bsigma}:=(\sigma_1,\ldots, \sigma_m)$  of Schubert conditions 
for $\Fl(\adot;n)$ satisfying 
$\ell(\sigma_1)+\cdots + \ell(\sigma_m) = \dim(\adot)$. 

Given a Schubert problem $\bsigma$ for $\Fl(\adot;n)$ and complete flags
$\Fdot^1,\dots, \Fdot^m$, the intersection 
 \begin{equation}\label{eq:SchubIntersect}
   X_{\sigma_1}\Fdot^1\cap\cdots\cap X_{\sigma_m}\Fdot^m\
 \end{equation}
is an \demph{instance} of $\bsigma$.
When the flags are in general position, this intersection is transverse and
zero-dimensional~\cite{Kleiman}, and it consists of all flags $\Edot\in\Fl(\adot;n)$
having position $\sigma_i$ relative to $\Fdot^i$, for each 
$i=1,\ldots, m$.
Such a flag $\Edot$ is a \demph{solution} to this instance of $\bsigma$.

The degree of a zero-dimensional intersection (\ref{eq:SchubIntersect}) is independent of 
the choice of the flags and we call this number $d(\bsigma)$ the
\demph{degree} of the Schubert problem $\bsigma$. 
When the intersection is transverse, the number of solutions to $\bsigma$ equals its degree. 

When the flags $\Fdot^1,\dots, \Fdot^m$ are real, the solutions to the Schubert
problem need not be real. 
The monotone secant conjecture posits a method to select the flags $\Fdot$ so that all
solutions are real, for a certain class of Schubert problems.

Let $\gamma\colon\R\to\R^n$ be a rational normal curve, which is any curve affinely equivalent
to the moment curve $t\mapsto(1,t,t^2,\ldots,t^{n-1})$. 
Consider this projectively, so that $\gamma$ is homeomorphic to
$\R\P^1$, which is a circle.
A flag $\Fdot$ is \demph{secant along an interval $I$} of $\gamma$ if every subspace in
the flag is spanned by its intersection with $I$. 
A list of flags $\Fdot^1,\dotsc,\Fdot^m$ secant to $\gamma$ is \demph{disjoint} if the
intervals of secancy are pairwise disjoint.
Disjoint flags are naturally ordered by the order in which their intervals of secancy lie
within $\R\P^1$.
We remark that this order is to be taken cyclically as in Figure~\ref{F:schematic}, and
with respect to one of the two orientations of $\R\P^1$.

A permutation $\sigma$ is \demph{Grassmannian} if it has a unique 
descent, for these, let $\delta(\sigma)$ be the descent.
A \demph{Grassmannian Schubert problem} is one that involves only Grassmannian Schubert
conditions. 
A list of disjoint secant flags $\Fdot^1,\dotsc,\Fdot^m$ is \demph{monotone} with respect to 
a Grassmannian Schubert problem $(\sigma_1,\dotsc,\sigma_m)$ if the function $\Fdot^i
\mapsto \delta(\sigma_i)$ is monotone with respect to one of the two orientations of
$\R\P^1$.
In other words, if
\[
   \delta(\sigma_i)\, <\, \delta(\sigma_j)\ 
    \Longrightarrow\ F^i < F^j\,,\qquad\mbox{for all } i,j\,,
\]
where $<$ is induced by one of the two cyclic orderings of $\R\P^1$.

\begin{monsecconj}\label{C:monosec}
 For any Grassmannian Schubert problem $(\sigma_1,\ldots, \sigma_m)$ on the flag manifold
 $\Fl(\adot;n)$ and any disjoint secant flags $\Fdot^1,\ldots,\Fdot^m$ that are monotone with
 respect to the Schubert problem, the intersection
\[
   X_{\sigma_1}\Fdot^1\cap X_{\sigma_2}\Fdot^2\cap
    \dotsb \cap X_{\sigma_m}\Fdot^m
\]
 is transverse with all points real.
\end{monsecconj}

Conjecture~\ref{C:first} is the monotone secant conjecture for a
Schubert problem on $\Fl(2,3;5)$ involving the Schubert conditions 
$\sigma:=13245$ and $\tau:=12435$, where we write permutations in one-line notation, so that
$\sigma(2)=3$ and $\tau(2)=2$.
Then $\delta(\sigma)=2$, $\delta(\tau)=3$, and $\ell(\sigma)=\ell(\tau)=1$,
so that $(\sigma,\sigma,\sigma,\sigma,\tau,\tau,\tau,\tau)$ is a
Schubert problem for $\Fl(2,3;5)$, as $\dim(\Fl(2,3;5))=8$.
We use exponential notation for repeated conditions, so that this Schubert problem is written
as \DeCo{$(\sigma^4,\tau^4)$}.
The corresponding Schubert varieties are
 \begin{eqnarray*}
   X_{\sigma}\Fdot\ =\ \{ \Edot\in\Fl(2,3;5) \mid \dim E_2\cap F_3 \geq 1 \}\,,\\
   X_{\tau}\Fdot\ =\ \{ \Edot\in \Fl(2,3;5) \mid \dim E_3 \cap F_2 \geq 1 \}\,,
 \end{eqnarray*}
that is, the set of flags $\Edot$ whose $2$-plane $E_2$ meets a fixed $3$-plane $F_3$
non-trivially,  and the set of $\Edot$ where $E_3$ meets a fixed $2$-plane $F_2$
non-trivially, respectively. 
Consequently, we write $X_\sigma F_3$ for $X_\sigma\Fdot$ and $X_\tau F_2$ for
$X_\tau\Fdot$. 
 
For $s,t,u,v,w\in\R$, let \DeCo{$F_3(s,t,u)$} be the linear span of $\gamma(s), \gamma(t)$, and
$\gamma(u)$ and \DeCo{$F_2(v,w)$} be the linear span of $\gamma(v)$ and $\gamma(w)$; these are a
secant 3-plane and a secant 2-plane to the rational normal curve, respectively. 
The condition $f(s,t,u;x)=0$ of Conjecture~\ref{C:first} is equivalent to the membership
$\Edot\in X_\sigma F_3(s,t,u)$.
Similarly, the condition $g(v,w;x)=0$ is equivalent to the membership $\Edot\in X_\tau F_2(v,w)$.
Lastly, the condition on the ordering of the points $s_i,t_i,u_i,v_j,w_j$ in
Conjecture~\ref{C:first} implies that the relevant subspaces $F_3(s_i,t_i,u_i)$ and
$F_2(v_j,w_j)$ for $i,j=1,\dotsc,4$ lie in disjoint secant flags that are monotone with respect
to this Schubert problem. 

Three conjectures that have driven progress in enumerative real algebraic
geometry are specializations of the monotone secant conjecture.
For the Grassmannian $\Gr(a;n)$, any list of disjoint secant flags
$\Fdot^1,\dotsc,\Fdot^m$ is monotone with respect to any Schubert problem
$(\sigma_1,\dotsc,\sigma_m)$, as all the conditions have the same descent.  
In this way, the monotone secant conjecture reduces to the secant conjecture.

\begin{sconj}\label{C:sec}
For any Schubert problem $(\sigma_1, \dots, \sigma_m)$ on any Grassmannian and any
disjoint secant flags $\Fdot^1,\ldots, \Fdot^m$, the intersection  
\[
  X_{\sigma_1}\Fdot^1\cap X_{\sigma_2}\Fdot^2\cap \cdots \cap X_{\sigma_m}\Fdot^m
\]
is transverse with all points real.
\end{sconj}

We studied this conjecture in a large-scale experiment whose results are described
in~\cite{FRSC_Sec}, solving 1,855,810,000 instances of 703 Schubert problems  on 13
different Grassmannians, verifying the secant conjecture in each of the 448,381,157
instances checked.
This took 1.065 teraHertz-years of computing.

The \demph{osculating flag $\Fdot(t)$} is the flag whose  $j$-dimensional
subspace $F_j(t)$ is the span of the first $j$ derivatives 
$\gamma(t),\gamma'(t),\ldots,\gamma^{(j-1)}(t)$ of $\gamma$ at $t$.  
This subspace $F_j(t)$ is the unique $j$-dimensional subspace having maximal 
order of contact, namely $j$, with $\gamma$ at $\gamma(t)$.
It follows that the limit of any family of flags whose intervals of secancy shrink to a point $\gamma(t)$
is this osculating flag $\Fdot(t)$.
In this way, the limit of the monotone secant conjecture, as the secant flags become
osculating flags, is a similar conjecture where we replace monotone secant flags by
monotone osculating flags.

\begin{monconj}\label{C:monotone}
 For any Schubert problem $(\sigma_1,\ldots, \sigma_m)$ on the flag manifold
 $\Fl(\adot;n)$ and any flags $\Fdot^1,\ldots, \Fdot^m$ osculating a rational normal
 curve $\gamma$ at real points that are monotone with respect to the Schubert problem, the
 intersection  
 \[
   X_{\sigma_1}\Fdot^1\cap X_{\sigma_2}\Fdot^2\cap \cdots\cap X_{\sigma_m}\Fdot^m
 \]
 is transverse with all points real.
\end{monconj}

Ruffo, et al.~\cite{RSSS} formulated and studied this conjecture, establishing special
cases and giving substantial experimental evidence in support of it.

The Shapiro conjecture is a specialization of the monotone secant conjecture that both
restricts to the Grassmannian and to osculating flags.
This was posed around 1995 by Boris 
Shapiro and Michael Shapiro and studied in~\cite{So_Shap}.
Proofs were given by 
Eremenko and Gabrielov for $\Gr(n{-}2;n)$~\cite{EG_02,EG11} using complex analysis and in
complete generality by Mukhin, Tarasov, and Varchenko~\cite{MTV_Annals,MTV_JAMS} using 
integrable systems and representation theory.

\begin{Shapiroconj}\label{ShapiroConj}
  For any Schubert problem $(\sigma_1, \ldots, \sigma_m)$ in a Grassmannian and any distinct
  real numbers $t_1,\ldots, t_m$, the intersection 
 \[
  X_{\sigma_1}\Fdot(t_1)\cap X_{\sigma_2}\Fdot(t_2)\cap\cdots\cap X_{\sigma_m}\Fdot(t_m)
 \]
 is transverse with all points real.
\end{Shapiroconj}

\section{Results}\label{Sec:results}

A consequence of the example discussed in Section~\ref{Sec:fourlines} is that
the secant conjecture (like the Shapiro conjecture before it) does not hold for flag
manifolds $\Fl(\adot;n)$ that are not Grassmannians. 
The monotonicity condition appears to correct this failure in both conjectures.
We give more details on the relation of the monotone conjecture to the monotone
secant conjecture and give a conjecture that interpolates between the two.
Then we discuss some of our data in an experiment that tested both conjectures. 

\subsection{The monotone conjecture is the limit of the monotone secant conjecture}\label{Sec:LimitMono}

The osculating plane $F_i(s)$ is the unique $i$-dimensional subspace having maximal order
of contact with the rational normal curve $\gamma$ at the point $\gamma(s)$, and therefore
it is a limit of secant planes.
We give a more precise statement of this fact.

\begin{proposition}\label{P:limit}
 Let $\{s_1^{(j)}\ldots, \ldots, s_i^{(j)}\}$ for $j=1,2,\ldots$ be a sequence of lists of
 $i$ distinct complex numbers with the property that for each $p=1,\ldots, i$, we have 
 \[
   \lim_{j\to \infty} s_p^{(j)}\ =\ s\,,
 \]
for some number $s$. Then,
 \[
   \lim_{j\to \infty} \Span\{ \gamma(s_1^{(j)}), \ldots, \gamma(s_i^{(j)}) \}\ =\ F_i(s)\,.
 \]
\end{proposition}

As explained in the previous section, by this proposition, the monotone secant conjecture 
implies the monotone conjecture.
This has a partial converse which follows from a standard limiting argument.

\begin{theorem}\label{T:limit}
 Let $(\sigma_1,\ldots, \sigma_m)$ be a Schubert problem on $\Fl(a;n)$ for which the
 monotone conjecture holds.
 Then, for any distinct real numbers that are monotone with respect to
 $(\sigma_1,\ldots,\sigma_m)$, there exists a number $\epsilon>0$ such that, if for each
 $i=1,\ldots,m$, $\Fdot^i$ is a flag secant to $\gamma$ along an interval of length $\epsilon$
 containing $t_i$, then the intersection 
 \[
   X_{\sigma_1}\Fdot^1\cap X_{\sigma_2}\Fdot^2\cap \cdots \cap X_{\sigma_m}\Fdot^m
 \]
is transverse with all points real.
\end{theorem}

%
\subsection{Generalized monotone secant conjecture}
We generalize the monotone secant conjecture, replacing secant flags by flags which are
spanned by osculating subspaces, as in~\cite[\S~3.3]{FRSC_Sec}.
By Proposition~\ref{P:limit}, such flags are intermediate between secant and osculating
flags, so this new conjecture interpolates between the monotone secant and monotone
conjectures. 
A \demph{generalized secant subspace} to the rational normal curve $\gamma$ is a subspace
that is spanned by subspaces osculating $\gamma$ at real points.
This notion includes secant subspaces as well as osculating subspaces, as a point of
$\gamma$ generates a one-dimensional subspace osculating $\gamma$.
A \demph{generalized secant flag} is one consisting of generalized secant subspaces.
A generalized secant flag is \demph{secant along an interval} $I$ if the osculating
subspaces spanning its subspaces osculate $\gamma$ at points of $I$.

\begin{conj}[Generalized monotone secant conjecture]
 For any Grassmannian Schubert problem $(\sigma_1,\ldots, \sigma_m)$ on the flag manifold
 $\Fl(\adot;n)$ and any disjoint generalized secant flags $\Fdot^1,\ldots,\Fdot^m$ that
 are monotone with respect to the Schubert problem, the intersection
\[
   X_{\sigma_1}\Fdot^1\cap X_{\sigma_2}\Fdot^2\cap
    \dotsb \cap X_{\sigma_m}\Fdot^m
\]
 is transverse with all points real.
\end{conj}

This conjecture contains the monotone and monotone secant conjectures as special cases, and
interpolates between the two.

\subsection{Experimental evidence for the monotone secant conjecture}\label{Sec:ExpEvidence}

While its relation to existing conjectures led to the monotone secant
conjecture, we believe the immense weight of empirical evidence is the strongest support
for it.
We conducted an experiment that tested 11,141,897,000 instances of 1300 Schubert problems
on 19 flag manifolds. 
Of these, 768,846,000 were instances of the monotone secant conjecture, which was
verified in every case tested.
We also tested 918,902,000 instances of the monotone conjecture for comparison. 
The remaining instances involved non-monotone evaluations of either
disjoint secant flags or osculating flags.
Our data consistently displayed a striking inner border, and a number of Schubert problems
exhibited lower bounds on their numbers of real solutions. 
This experiment used 1.9 teraHertz-years of computing.

Table \ref{T:X1^4Y1^4=12} shows the data we obtained for the Schubert problem
$(\Blue{\sigma^4},\G{\tau^4})$ with 12 solutions on the flag manifold $\Fl(2,3;5)$
  introduced in Conjecture~\ref{C:first}.  
\begin{table}[htb]
\noindent{\small 

\noindent\begin{tabular}{r|r||r|r|r|r|r|r|c||r|}
\multicolumn{10}{c}{Real Solutions}\\
\cline{2-10}
\multirow{10}{*}{\begin{sideways}Necklace\end{sideways}} &
 & 0&2&4&6&8&10&12 &Total\\\cline{2-10}\noalign{\smallskip}\cline{2-10}
&\Blue{2222}\G{3333} &&&&&&&1000000 & 1000000\\\cline{2-10}
&\Blue{222}\G{33}\Blue{2}\G{33} &&&31&205380&545269&1245017&1004303& 3000000\\\cline{2-10}
& \Blue{2}\Blue{2}\G{3}\Blue{2}\Blue{2}\G{3}\G{3}\G{3} &&&196&403485&1071579&967226&557514& 3000000\\\cline{2-10}
& \Blue{2}\Blue{2}\G{3}\G{3}\Blue{2}\Blue{2}\G{3}\G{3} &&&391&801525&1200700&651183&346201& 3000000\\\cline{2-10}
& \Blue{2}\Blue{2}\G{3}\Blue{2}\G{3}\G{3}\Blue{2}\G{3} &&1025&70009&300121&938430&1123770&566645& 3000000\\\cline{2-10}
& \Blue{2}\Blue{2}\G{3}\Blue{2}\G{3}\Blue{2}\G{3}\G{3} &&33488&950203&1256341&560205&164020&35743& 3000000\\\cline{2-10}
& \Blue{2}\Blue{2}\Blue{2}\G{3}\Blue{2}\G{3}\G{3}\G{3} &&93232&284316&460016&1010425&750171&401840&3000000\\\cline{2-10}
& \Blue{2}\G{3}\Blue{2}\G{3}\Blue{2}\G{3}\Blue{2}\G{3} &885953&854550&830122&298843&104574&23741&2217& 3000000\\\cline{2-10}\noalign{\smallskip}\cline{2-10}
&Total &885953&982295&2135268&3725711&5431182&4925128&5914463 & 24000000 \\\cline{2-10}
\end{tabular}

}\vspace{5pt}
\caption{Necklaces vs. real solutions for $(\Blue{\sigma^4},\G{\tau^4})$ in $\Fl(2,3;5)$.}\vspace{-20pt}
\label{T:X1^4Y1^4=12}
\end{table}
We computed 24,000,000 instances of this problem, all involving flags that were secant to
the rational normal curve along disjoint intervals. 
This took 17.618  gigaHertz-years. 
The columns are indexed by even integers from 0 to 12, indicating the number of
real solutions. 
The rows are indexed by \demph{necklaces}, which are sequences
$\delta(\sigma_1),\ldots,\delta(\sigma_m)$, where $\delta(\sigma_i)$ denotes the
descent of the Grassmannian permutation $\sigma_i$, as described in
Section~\ref{Sec:background}.  
In the table a \Blue{2} represents the condition on the two-plane $E_2$
given by the permutation $\Blue{\sigma} =\Blue{13245}$, and a \G{3} represents
the condition on $E_3$ given by the permutation $\G{\tau} =\G{12435}$.  

In Table \ref{T:X1^4Y1^4=12}, the first row labeled with $\Blue{2222}\G{3333}$ represents
tests of the monotone secant conjecture, since the only entries are in the column for 12
real solutions, the monotone secant conjecture was verified in 3,000,000 instances. 
This is the only row with only real solutions. 

Compare this to the 16,000,000 instances of the same Schubert problem, but with osculating
flags, which we present in Table \ref{T:2X1^4Y1^4=12}. 
This computation took 85.203
gigahertz-days. 
\begin{table}[htb]
\noindent{\small 

\noindent\begin{tabular}{r|r||r|r|r|r|r|r|c||r|}
\multicolumn{10}{c}{Real Solutions}\\
\cline{2-10}
\multirow{10}{*}{\begin{sideways}Necklace\end{sideways}} &
 & 0&2&4&6&8&10&12 &Total\\\cline{2-10}\noalign{\smallskip}\cline{2-10}
&\Blue{2222}\G{3333} 
&&&&&&&2000000 & 2000000\\\cline{2-10}
&\Blue{222}\G{33}\Blue{2}\G{33}
 &&&1041&246876&581972&582865&587246& 2000000\\\cline{2-10}
& \Blue{2}\Blue{2}\G{3}\Blue{2}\Blue{2}\G{3}\G{3}\G{3}
 &&&1480&263981&621920&584508&528111& 2000000\\\cline{2-10}
& \Blue{2}\Blue{2}\G{3}\G{3}\Blue{2}\Blue{2}\G{3}\G{3}
 &&&8882&217100&861124&503562&409332& 2000000\\\cline{2-10}
& \Blue{2}\Blue{2}\G{3}\Blue{2}\G{3}\G{3}\Blue{2}\G{3}
&&\begin{picture}(1,1)\put(-30,-2){\Yellow{\rule{35pt}{9pt}}}\end{picture}
  &120195&402799&665766&549653&261587& 2000000\\\cline{2-10}
& \Blue{2}\Blue{2}\G{3}\Blue{2}\G{3}\Blue{2}\G{3}\G{3} 
&&7329&255114&431074&664551&420699&221233& 2000000\\\cline{2-10}
& \Blue{2}\Blue{2}\Blue{2}\G{3}\Blue{2}\G{3}\G{3}\G{3} 
&&25552&137116&227922&415553&424582&769275& 2000000\\\cline{2-10}
& \Blue{2}\G{3}\Blue{2}\G{3}\Blue{2}\G{3}\Blue{2}\G{3} 
& 49197&125725&557851&395992&516675&244212&110348& 2000000\\\cline{2-10}\noalign{\smallskip}\cline{2-10}
&Total & 49197&158606&1081679&2185744&4327561&3310081&4887132& 16000000 \\\cline{2-10}
\end{tabular}

}\vspace{5pt}
\caption{Necklaces vs. real solutions for  $(\Blue{\sigma^4},\G{\tau^4})$ in $\Fl(2,3;5)$.}\vspace{-20pt}
\label{T:2X1^4Y1^4=12}
\end{table}
Both tables are similar with nearly identical ``inner borders'', except for the shaded box
in Table~\ref{T:2X1^4Y1^4=12}.
In fact, by a standard argument similar to that which implied Theorem~\ref{T:limit}, we
may conclude that every number of real solutions to a Schubert problem observed for a
given necklace with osculating flags also occurs for that Schubert problem and necklace 
with secant flags, where the points of secancy are sufficiently clustered.
That is, the support of a table for the monotone conjecture will be a subset of the
support of the corresponding table for the monotone secant conjecture.
There are some Schubert problems for which we did not observe this containment; the reason
for this is that we aparently did not compute an instance with secant flags whose points of secancy
were sufficiently clustered.

%
\subsection{Lower bounds and inner borders}\label{Sec:inner} 
 
The most enigmatic phenomenon that we observe in our data is the presence of an
``inner border" for many geometric problems, as we have pointed out in example of
Table~\ref{T:X1^4Y1^4=12}.
That is, for some necklaces (besides the monotone ones), there appears to be a lower bound
on the number of real solutions.
We do not understand this phenomenon, even conjecturally.
Our software that displays the tables is designed to highlight this feature of our data.

Another common phenomenon is that for many problems, there are always at least some
real solutions, for any necklace.
(Note that the last rows of Tables~\ref{T:X1^4Y1^4=12} and~\ref{T:2X1^4Y1^4=12} had
instances with no real solutions).
Table~\ref{T14:W3X4=21_MSC} displays an example of this for 
a Schubert problem $(\Purple{\sigma^3},\Blue{\tau^5})$ on $\Fl(2,3;6)$ 
\begin{table}[htb]
{\small
\[
  \begin{tabular}{r|r||r|r|r|r|r|r|r|r|r|r|r|r|}
\multicolumn{14}{c}{Real Solutions}\\
\cline{2-14}  
\multirow{7}{*}{\begin{sideways}Necklace\end{sideways}} &
    &1&3&5&7&9&11&13&15&17&19&21&Total\\\cline{2-14}\noalign{\smallskip}\cline{2-14}
&   \W\W\W\X\X\X\X\X&&&&&&&&&&&80000&80000\\\cline{2-14}
 &  \W\W\X\W\X\X\X\X&&&&&&&921&16549&26267&14475&21788&80000\\\cline{2-14}
  & \W\W\X\X\W\X\X\X&&&&&&39&1208&24559&39013&13947&1234&80000\\\cline{2-14}
   &   \W\X\W\X\X\W\X\X&&&&&&612&9544&43256&23583&2927&78&80000\\\cline{2-14}
  & \W\X\W\X\W\X\X\X&&&&&&3244&19887&31931&13688&3632&7618&80000\\\cline{2-14}
   &Total&&&&&&3895&31560&116295&102551&34981&110718&400000\\\cline{2-14}
  \end{tabular}
\]	 
}\caption{Enumerative Problem $W^3X^5= 21$ on $\Fl(2, 3 ; 6)$}\vspace{-20pt}
\label{T14:W3X4=21_MSC}
\end{table}
with 21 solutions, where
$\Purple{\sigma}:=142356$ has $\delta(\Purple{\sigma})=2$ and 
$\ell(\Purple{\sigma})=2$ and $\Blue{\tau}:=124356$ has $\delta(\Blue{\tau})=3$ 
and $\ell(\Blue{\tau})=1$.
Very prominently, it appears that at least 11 of the solutions will always be real.
This computation took 7.67 gigaHertz-years.

Such lower bounds and inner borders were observed when studying the
monotone conjecture~\cite{RSSS}.
Eremenko and Gabrielov established lower bounds for the Wronski map~\cite{EG01} in
Schubert calculus for the Grassmannian, more recently Azar and
Gabrielov~\cite{AzarGab} established lower bounds for some instances of the monotone
conjecture which were observed in~\cite{RSSS}.

%
\section{Methods}\label{Sec:method}

Our experimentation was possible as instances of Schubert problems are simple to model on a
computer. 
The procedure we use may be semi-automated and run on supercomputers. 
We will not describe how this automation is done, for that is the subject of the
paper~\cite{Exp-FRSC}; instead, we explain here the computations we performed.

For a Schubert condition $\sigma$ on a flag variety $\Fl(\adot;n)$, let \DeCo{$j(\sigma)$} be
the dimension of the largest subspace in a flag $\Fdot$ that is needed to define
$X_{\sigma}\Fdot$. 
For example, we have seen that $j(13245)=3$ and $j(12435)=2$.

To compute an instance of a Schubert problem $(\sigma_1, \ldots, \sigma_m)$ corresponding to a
necklace $\nu$, we first select $N:=N(\sigma_1)+\dotsc+N(\sigma_m)$ real numbers and then group
them into disjoint subsets $s^{(1)},\dotsc,s^{(m)}$ where $s^{(i)}$ consists of $N(\sigma_i)$
consecutive numbers.
Furthermore, the relative ordering of the subsets corresponds to the necklace $\nu$.
Having done this, each subset $s^{(i)}$ defines a secant (partial) flag $\Fdot(s^{(i)})$.
We use these flags to formulate the instance of the Schubert problem
\[
   X_{\sigma_1}\Fdot(s^{(1)})\,\cap\,
   X_{\sigma_2}\Fdot(s^{(2)})\,\cap\,\dotsb\,\cap\,
   X_{\sigma_m}\Fdot(s^{(m)})
\]
as a system of polynomials in $\dim(\adot)$ local coordinates for $\Fl(\adot;n)$, whose common
zeroes represent the solutions to this instance of the Schubert problem.
This was illustrated in the Introduction when Conjecture \ref{C:first} was presented. 
See~\cite{Fu97,RSSS, So_Shap} for details. 

Given this system of polynomials, we use Gr\"obner bases to compute a polynomial in one variable of minimal degree
in the ideal of these equations. 
This univariate polynomial is called an \demph{eliminant}. 
If the eliminant is square-free and has degree equal to the expected number of complex
solutions (this is easily verified) then the number of real
roots of the eliminant equals the number of real solutions to the Schubert problem.  
This follows from the form of a lexicographic Gr\"obner basis for the ideal, as described by the Shape
Lemma~\cite{BMMT}. 
This is given in more detail in \S 2.2 of~\cite{IHP}.
We determine the number of real solutions to the Schubert problem by computing the number of real roots of the
eliminant. 
For this, we use MAPLE's {\tt realroot} command, which uses symbolic methods based on Sturm sequences to determine
the number of real roots of a univariate polynomial.
If the software is reliably implemented, which we believe, then this computation provides a
proof that the given instance has the computed number of real solutions to the original
Schubert problem.  

In our computations, for a given Schubert problem, we first make a choice of $N$ real numbers,
and then use these same $N$ numbers for all necklaces for that problem.
Then we make another choice, and so on, ultimately making thousands to millions of such
choices.

For each Schubert problem we studied, we not only tested many instances of the monotone secant
conjecture, but also of the monotone conjecture, comparing the two as we did for
the Schubert problem $(\Blue{\sigma^4},\G{\tau^4})$ in $\Fl(2,3;5)$, where $\Blue{\sigma}=13245$ and
$\G{\tau}=12435$.
To compute instances of the monotone conjecture, we choose real numbers $s_1,\dotsc,s_m$ and
use osculating flags $\Fdot(s_1),\dotsc,\Fdot(s_m)$.
This is also described in~\cite[\S~5]{RSSS}.

For many Schubert problems, it was infeasible to compute instances of the monotone secant
conjecture, and we instead computed instances of the generalized monotone secant conjecture.
For these, one of the flags was the flag $\Fdot(\infty)$ osculating the rational normal curve
at infinity.
Then we used local coordinates for $X_{\sigma_1}\Fdot(\infty)$, in place of the local
coordinates for $\Fl(\adot;n)$; this uses $\ell(\sigma_i)$ fewer local coordinates.

For some Schubert problems we wanted to study, there were several hundred to many thousands of
necklaces, and for these we uniformly chose a much smaller set of necklaces to compute, which
we called coarse necklaces.
In our on-line tables, we encoded these choices in a variable called computation type.
Computation types 4 and 7 were for instances of the monotone conjecture, 5 and 8 for the
monotone secant conjecture, and 6 and 9 for the generalized monotone secant conjecture.
In each of these, the first number indicates that we used all necklaces, while the second we
used coarse necklaces.



\providecommand{\MR}{\relax\ifhmode\unskip\space\fi MR }
\providecommand{\MRhref}[2]{%
  \href{http://www.ams.org/mathscinet-getitem?mr=#1}{#2}
}
\providecommand{\href}[2]{#2}


\end{document}